\documentclass[12pt,reqno]{article}

\usepackage[usenames]{color}
\usepackage[colorlinks=true,
linkcolor=webgreen,
filecolor=webbrown,
citecolor=webgreen]{hyperref}
\definecolor{webgreen}{rgb}{0,.5,0}
\definecolor{webbrown}{rgb}{.6,0,0}

\usepackage{color}
\usepackage{epsfig}
\usepackage{amsmath,amssymb}
\usepackage{amsfonts}
\usepackage{appendix}

\newcommand{\seqnum}[1]{\href{http://www.research.att.com/cgi-bin/access.cgi/as/~njas/sequences/eisA.cgi?Anum=#1}{\underline{#1}}}

\newenvironment{packed_enumerate}{
\setlength{\parsep}{0pt}
\setlength{\parskip}{0pt}
\begin{enumerate}
  \setlength{\itemsep}{1pt}
  \setlength{\parsep}{0pt}
  \setlength{\parskip}{0pt}
}{\end{enumerate}}

\begin{document}

\begin{center}
\epsfxsize=4in
\leavevmode\epsffile{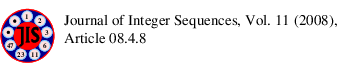}
\end{center}

\begin{center}
\vskip 1cm{\LARGE\bf 
Solving Triangular Peg Solitaire}
\vskip 1cm
\large
George I. Bell\\
Tech-X Corporation\\
5621 Arapahoe Ave, Suite A\\
Boulder, CO 80303 \\
USA\\
\href{mailto:gibell@comcast.net}{\tt gibell@comcast.net} \\
\end{center}

\vskip .2 in
\begin{abstract}
We consider the one-person game of peg solitaire
on a triangular board of arbitrary size.
The basic game begins from a full board with one peg missing
and finishes with one peg at a specified board location.
We develop necessary and sufficient conditions for this
game to be solvable.
For all solvable problems,
we give an explicit solution algorithm.
On the 15-hole board, we compare three simple solution strategies.
We also consider the problem of finding solutions that minimize the
number of moves (where a move is one or more
consecutive jumps by the same peg),
and find the shortest solution to the basic game
on all triangular boards with up to 55 holes (10 holes on a side).
\end{abstract}

\newtheorem{theorem}{Theorem}[section]
\newtheorem{proposition}{Proposition}[section]
\newtheorem{corollary}{Corollary}[section]
\newtheorem{lemma}{Lemma}[section]

\section{Introduction} 
\begin{figure}[htb]
\centering
\epsfig{file=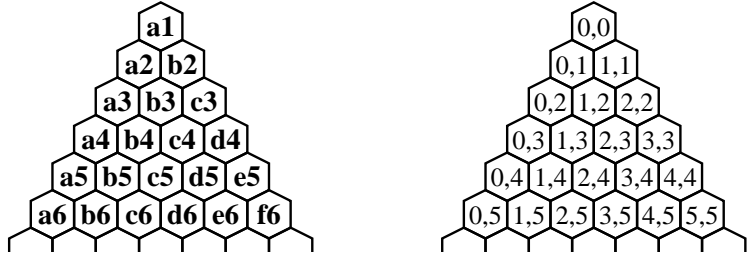}
\caption{The board $T_n$ with two types of hole coordinates (alphanumeric and skew).}
\label{fig1}
\end{figure}

For many years, $\mbox{Cracker Barrel}\textsuperscript{\textregistered}$~restaurants
have popularized peg solitaire played on a triangular board with
15~holes.
Many patrons have puzzled over this game, often called ``an IQ test",
which is surprisingly difficult given its small size and simple rules.
Often people resort to a computer program to solve this puzzle,
and it is a popular assignment in computer science classes \cite[p. 132]{CC}.
In this article we consider peg solitaire on a triangular board
with $n$ holes on each side.
This board will be referred to as $T_n$ and can be
conveniently presented on an array of hexagons (Figure~\ref{fig1}).
The board $T_n$ has $T(n)=n(n+1)/2$ holes,
where $T(n)$ is the $n$th triangular number.
The $\mbox{Cracker Barrel}\textsuperscript{\textregistered}$ board is $T_5$.

We will use two different notations to identify the holes in the board.
The notation in Figure~\ref{fig1}a is useful for quick hole identification
and for describing solutions.
The ``skew Cartesian coordinate" notation \cite{Duncan} in Figure~\ref{fig1}b is useful
for the theory of the game,
as well as inside computer programs.
It is also particularly easy to perform reflections and rotations
of the board in this coordinate system (see Appendix~\ref{rotref}).
Note that these two coordinate notations are closely related\footnote{To convert from skew coordinates to the alphanumeric notation,
map the $x$-skew coordinate to the alphabet
($0\rightarrow a$, $1\rightarrow b$, $\ldots$)
and concatenate with the $y$-skew coordinate plus one.}.

The game begins with a peg (or marble) at every hole except one,
called the \textbf{starting vacancy}.
The player then \textbf{jumps} one peg over another into an empty hole on
the board, removing the peg that was jumped over.
The game ends when no jump is possible,
and the goal is to finish at a one peg position.
If the starting vacancy and the ending hole happen to be the same hole,
then we call this a
\textbf{complement problem}\footnote{This type of problem has also been called a ``reversal" \cite{WinningWays}.}.
For example, a popular $T_5$ puzzle is to start with one peg missing in the top corner,
and try to finish with one peg in the same corner, the a1-complement problem.
We will use the term \textbf{move} for one or more consecutive jumps by the same peg.
To denote a jump, we will list the starting and ending coordinates
separated by a dash, i.e., a1-a3.
When the same peg makes two or more jumps in a single move,
instead of listing each jump separately (a1-a3, a3-c3, c3-a1)
we will combine them by writing a1-a3-c3-a1.

This puzzle is a variant of \textbf{square lattice} solitaire,
generally played on a 33-hole cross-shaped board \cite{BellFr}.
Both puzzles have the same jumping rules,
with the 33-hole board formed from a square lattice of holes,
while triangular solitaire is played
on a triangular (or hexagonal) lattice of holes.
Square lattice solitaire has a 300 year history,
but the origins of triangular solitaire are more obscure.
Triangular solitaire was popularized by a 1966
Martin Gardner column \cite{Gardner},
where he considered the game played using
a triangular array of pennies on a table.
However, an 1891 patent \cite{patent} indicates that
triangular solitaire is quite a bit
older\footnote{This patent is for a 16-hole board on a triangular lattice,
but it is not a triangular board.}.
Hentzel \cite{Hentzel1} published the first mathematical analysis of
the game in 1973.

\section{Theory of the Game}

In square lattice solitaire,
the set of finishing holes is restricted by the so-called
\textbf{rule of three} \cite{WinningWays}, which states
that the $x$-coordinates of possible finishing holes
differ by a multiple of 3 (and similarly for the $y$-coordinates).
The analog of this theory for triangular solitaire
was given by Hentzel \cite{Hentzel1} in 1973.
Generally, this theory is developed using an elegant
group-theoretic argument \cite{Hentzel1, WinningWays, Duncan}.
We will use a simpler parity argument to prove our results.

\subsection{The Four Position Classes}

\begin{theorem}\label{th1}
On the triangular board $T_n$ with $n\ge 4$,
beginning from a vacancy with skew-coordinates $(x_s,y_s)$,
the following conditions are equivalent:\newline
\indent A) The board is not solvable to one peg;\newline
\indent B) $n\equiv 1 \pmod 3$ and $x_s+y_s \equiv 0 \pmod 3$.
\end{theorem}

\begin{figure}[htb]
\centering
\epsfig{file=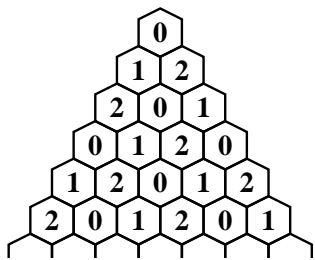}
\caption{The labeling of holes for the parity argument.
The hole with skew coordinates $(x,y)$ is labeled $(x+y) \bmod 3$.}
\label{fig2}
\end{figure}

{\it Proof} ($B\Rightarrow A$):
Label the holes in the board with the pattern in Figure~\ref{fig2},
where a hole with skew coordinates $(x,y)$ is labeled $(x+y) \bmod 3$.
This labeling pattern was chosen because every
jump involves exactly one hole of each of the three labels.
Let $c_i$ be the number of pegs in the holes labeled $i$.
After a jump is executed,
two of the three $c_i$ decrease by 1,
while the other increases by 1.
Therefore, if we add any pair among $\{c_0, c_1, c_2\}$,
the parity (even or odd) of this sum cannot change as the game is played.
We can represent the three invariant parities as the binary 3-vector
$(c_1+c_2,c_0+c_2,c_0+c_1)$,
where each component is taken modulo 2.

We can partition the set of all possible board positions into
four equivalence classes, called \textbf{position classes}.
During a solitaire game the board position remains in the same position class.
Moreover, there is a simple board position which can be chosen
as a representative of each position class.
The four position classes are defined by the
values of $(c_1+c_2,c_0+c_2,c_0+c_1)$, namely
$\mbox{EMPTY}=(0,0,0)$, $\mbox{PEG}_0=(0,1,1)$,
$\mbox{PEG}_1=(1,0,1)$ and $\mbox{PEG}_2=(1,1,0)$.
Other alternatives, such as $(1,1,1)$, can never occur,
because it is impossible to have all three sums odd,
or exactly one sum odd.

The empty board lies in the position class $\mbox{EMPTY}$,
while any board position with a single peg at a hole
labeled $i$ lies in class $\mbox{PEG}_i$
(in this sense we can refer to the label $i$ as
the \textit{position class of the hole}).
The position class $\mbox{EMPTY}$ is the only position class
which has no representative with a single peg.
Therefore, \textit{no board position in class} $\mbox{EMPTY}$
\textit{can be reduced to a single peg} using peg solitaire jumps.
So when is a board position with one peg missing in the
class $\mbox{EMPTY}$?
If $n \not\equiv 1 \pmod 3$ then the full board is in class $\mbox{EMPTY}$,
so if we remove a peg labeled $i$ in Figure~~\ref{fig2},
we are in class $\mbox{PEG}_i$, and can only finish with one peg
at a hole labeled $i$.

Triangular boards $T_n$ with $n\equiv 1 \pmod 3$
are the only triangular boards with a central hole,
and are also characterized by having a total number of holes not divisible by 3.
If $n\equiv 1 \pmod 3$ then the full board is in class $\mbox{PEG}_0$.
If we remove any peg labeled $0$,
the board is in position class $\mbox{EMPTY}$,
and cannot be solved to a single peg.
Any corner vacancy, as well as the central vacancy,
is always in position class $\mbox{EMPTY}$.
Thus,
\textit{no central vacancy problem is solvable on any triangular board},
a conclusion also reached by Beasley \cite[p. 231]{Beasley}.
It is easy to check that if we remove a peg labeled $1$ ($2$),
we are in class $\mbox{PEG}_2$ ($\mbox{PEG}_1$).
Thus if we remove a peg at $2$,
we can only finish with one peg at a hole labeled $1$,
and if we remove a peg at $1$,
we can only finish with one peg at a hole labeled $2$.

This proves $B\Rightarrow A$ in Theorem~\ref{th1}.
To complete the proof of Theorem~\ref{th1},
it suffices to show that any $T_n$ with $n\ge 4$ satisfying
$n\not\equiv 1 \pmod 3$ or
$n\equiv 1 \pmod 3$ and $x_s+y_s\not\equiv 0 \pmod 3$
can be solved down to one peg ($\sim B\Rightarrow \sim A$).
This part of the proof will be completed in Section~\ref{finish_th1}.

\vskip 12pt

Any starting vacancy which does not satisfy condition
B will be called a \textbf{feasible} starting vacancy,
because it can potentially be solved down to one peg.
We could also select a particular finishing hole and consider pairs of
starting and finishing holes that meet the above parity requirements,
called a \textbf{feasible pair}.
The following theorem tells us when
it is possible to play between a feasible pair of holes.

\begin{theorem}\label{th2}
Consider the triangular board $T_n$ with
starting vacancy $(x_s,y_s)$ and finishing hole $(x_f,y_f)$.
Then the following is a necessary condition for this
problem to be solvable:\newline
\indent 1) if $n\equiv 1 \pmod 3$, then $x_s+y_s\not\equiv 0 \pmod 3$ and $x_s+y_s+x_f+y_f\equiv 0 \pmod 3$, or\newline
\indent 2) if $n\not\equiv 1 \pmod 3$, then $x_s+y_s\equiv x_f+y_f \pmod 3$.\newline
In addition, for $n\ge 6$, the above condition is also sufficient.
\end{theorem}

{\it Proof}: That the condition is necessary is
a restatement of the parity arguments just presented.
For example, the condition $x_s+y_s\equiv x_f+y_f \pmod 3$ specifies that
the starting and finishing board positions must be in
the same position class.
For an alternative proof using an algebraic argument, see Duncan and Hayes \cite{Duncan}.
The sufficient part of the proof must show that when $n\ge 6$ any feasible pair
is in fact solvable, and this will be given in Section~\ref{finish_th2}.

\vskip 12pt

For a given board size $n$, we need some way of accounting for all
possible feasible pairs
(not duplicating pairs equivalent by reflection and/or rotation of the board).
A useful fact is that for $n\not\equiv 1 \pmod 3$,
we can cover all possible cases by considering
\textit{only starting vacancies and finishing holes in position class}
$\mbox{PEG}_0$.
The reason why this works is because each of the three corners
is in a different position class.
If $n\equiv 1 \pmod 3$, then we can take all starting vacancies in
position class $\mbox{PEG}_2$ (or $\mbox{PEG}_1$).

We conclude this section with a theorem about the symmetry
of board positions that can appear during solutions.

\begin{theorem}\label{th3}
Consider the triangular board $T_n$ with $n\not\equiv 1 \pmod 3$.
A solution to any problem finishing with one peg
cannot pass through a board position with $120^\circ$ rotational symmetry.
\end{theorem}

{\it Proof}: On $T_n$ with $n\not\equiv 1 \pmod 3$ it is easy to see
that any board position with $120^\circ$ rotational symmetry must
have $c_0=c_1=c_2$, so is in position class $\mbox{EMPTY}$.
Consequently, no rotationally symmetric board position can be
reduced to one peg.

\begin{corollary}\label{co1}
On any triangular board $T_n$, the solution to a complement
problem can never pass through a board position with $120^\circ$ rotational symmetry.
\end{corollary}

\subsection{Solving the Triangular Board \texorpdfstring{$T_4$}{T4}}

The 10-hole board $T_4$ is the smallest triangular board on which a
problem beginning with one peg missing is solvable to one peg.
This board falls in the first category of Theorem~\ref{th2}, and
this theorem gives us three geometrically distinct problems
that are potentially solvable:
beginning at a2 and finishing at b2, a3, or c4.

One interesting property of this board is that there is no
way to move a peg to the center b3.
Consequently, a solution must include exactly one jump over b3.
All other jumps originate or end at corners, and it is not
hard to see that two corners must have two jumps leaving them and one into them,
while the third corner has one jump leaving it.
This accounts for all eight jumps in a solution.

\begin{figure}[htb]
\centering
\epsfig{file=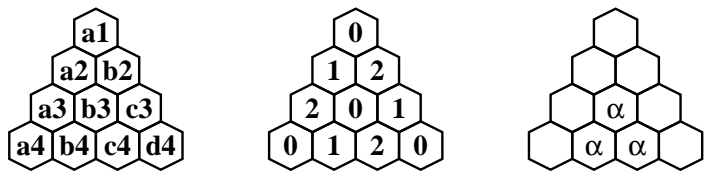}
\caption{$T_4$ notation (a), position class (b), and the parity count $\alpha$ (c).}
\label{fig3}
\end{figure}

Figure~\ref{fig3}c shows a useful parity count on this board.
The parity count $\alpha$ is the parity (even or odd) of the
number of pegs in the holes marked $\alpha$.
The only way to change this parity is using the jump a3-c3, or c3-a3.
There are two other similar parity counts obtained by rotating
the board.
It is easy to check that any solution from the
a2 vacancy to a3 must change \textit{all three}
parity counts.
This is impossible since the peg at b3 can only be removed
once---any solution changes exactly one of the three parities.

The other two problems do involve changing one of the parities,
and this implies that a solution to the problem
from a2 to b2 (if it exists) must contain the jump a3-c3 or c3-a3,
and a solution to the problem from a2 to c4 must contain a2-c4 or c4-a2.
The first problem is solvable (Appendix~\ref{Tn_sols}) but the second is not.
I have not found a simple argument showing
that the second problem is unsolvable.
However, it is easy to verify by calculating
the game tree by hand\footnote{Only 27 board positions can be reached starting with a2 vacant, assuming b3 must be cleared by a2-c4 or c4-a2.}
or using a computer program \cite{Miller}.

\subsection{Solving the Cracker Barrel Board \texorpdfstring{$T_5$}{T5}}

The $T_5$, or $\mbox{Cracker Barrel}\textsuperscript{\textregistered}$ board,
is one of the most interesting triangular boards.
It is (too) easy to write a short program to find solutions on this board.
Consider, for example, the problem starting and finishing in a corner, the a1-complement.
A program can find all solutions in a fraction of a second,
one of which is given in Appendix~\ref{Tn_sols}.
This solution has been converted\footnote{The holes are numbered sequentially starting with $1$,
each jump is specified by an ordered pair of the beginning and ending holes.
The solution is then a sequence of $13$ such ordered pairs, or 26 integers.}
to integer sequence \seqnum{A120422} in OEIS \cite{OEIS}.
A program can count that there are $6,816$ distinct solutions to the a1-complement.
A second solution is given by the sequence of jumps:
(a3-a1, c5-a3, a5-c5, d5-b5, c3-c5, b5-d5, e5-c5, a1-c3,
d4-b2, a4-a2,
b2-b4, c5-a3-a1).
In fact, \textit{any solution} to the a1-complement is either a permutation
of the jumps of one of these two solutions,
or its reflection about the $y$-axis is.
In this sense there are only \textit{two} solutions to the a1-complement.

\vskip 12pt
\begin{figure}[htb]
\centering
\epsfig{file=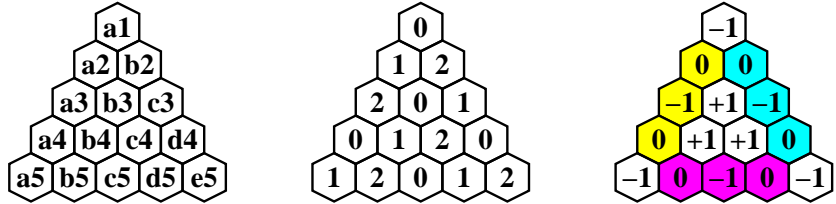}
\caption{$T_5$ notation (a), position class (b), and the SAX count (c).}
\label{fig4}
\end{figure}

Most people lose interest in this puzzle once they have found a solution.
However, additional theory yields many useful insights.
A pagoda function \cite{WinningWays, Beasley} is a real-valued function
of the board position that cannot increase as the game is played.
In 1986, Hentzel and Hentzel \cite{Hentzel2} discovered a powerful
pagoda function on this board.
Suppose we consider the hole-weighting shown in Figure~\ref{fig4}c,
and sum up the weights where a peg is present.
This function is not a pagoda function because a jump along the edge,
over one of the holes marked ``$-1$" can increase it.
However we can remedy this by adding ``$+1$" to our function for each
3-hole colored (or shaded) edge region which contains two or three pegs.
It is not hard to show that this function can never increase,
no matter what the starting board position is
and what jump is executed \cite{Hentzel2}.

As defined\footnote{Actually, Hentzel and Hentzel \cite{Hentzel2}
define their SAX count to be the negative of that defined here,
and prove that it can never decrease.}
by Hentzel and Hentzel \cite{Hentzel2},
we compute the SAX count as $S+A-X$ where:\newline
\indent $\bullet$ $S$ is the number of colored edge regions with two or more pegs
($0\le S \le 3$).\newline
\indent $\bullet$ $A$ is the number of pegs occupying holes labeled ``$+1$" in
Figure~\ref{fig4}c.\newline
\indent $\bullet$ $X$ is the number of pegs occupying holes labeled ``$-1$" in
Figure~\ref{fig4}c.\newline
Note that if the entire board is filled, the SAX count is $3+3-6=0$.
If a board position has only a single peg,
then the SAX count is simply the value of
that hole in Figure~\ref{fig4}c.
For any board position $B$, if we take the \textbf{complement} of $B$
(where every peg is replaced by a hole and vice versa), then
the SAX count of the complemented position is $-SAX(B)$.

Table~\ref{table1} shows all 17
distinct feasible starting and ending pairs on $T_5$.
The \textbf{slack} is the difference between the SAX count of the starting board
position and the SAX count of the ending board position.
In the case where we begin at a corner (say a1),
the starting SAX count is $+1$,
but the first jump must be a3-a1 (or the symmetric c3-a1),
and the SAX count is zero after this jump is made.
We define the \textbf{effective slack} as the difference
between the starting and final SAX count
when the effect of these forced jumps at the start or finish
is taken into account.
The effective slack is one less than the slack
when the game begins at a corner,
and one less when it ends at a corner.

\begin{table}[htb]
\begin{center} 
\begin{tabular}{| c | c | c | c || c | c | c | c |}
\hline
& Finish & Effective & & & Finish & Effective & \\
Vacate & At & Slack & Solvable? & Vacate & At & Slack & Solvable? \\
\hline
\hline
c5 & c5 & 2 & Yes & a4 & a1 & 0 & Yes\\ \hline
a1 & c5 & 1 & Yes & a4 & a4 & 0 & Yes\\ \hline
c5 & a1 & 1 & Yes & a4 & d4 & 0 & Yes\\ \hline
c5 & a4 & 1 & Yes & a1 & b3 & $-1$ & No\\ \hline
a4 & c5 & 1 & Yes & b3 & a1 & $-1$ & No\\ \hline
a1 & a1 & 0 & Yes & b3 & a4 & $-1$ & No\\ \hline
a1 & a4 & 0 & Yes & a4 & b3 & $-1$ & No\\ \hline
b3 & c5 & 0 & Yes & b3 & b3 & $-2$ & No\\ \hline
c5 & b3 & 0 & Yes & & & & \\ \hline
\end{tabular}
\caption{The 17 distinct feasible pairs on $T_5$ starting with
one peg missing and finishing with one peg.  Only 12 are solvable.
Note that all problems are in position class $\mbox{PEG}_0$.} 
\label{table1}
\end{center} 
\end{table}

The fact that the SAX count cannot increase proves
that any game with negative
effective slack cannot be solved.
This is 5 out of 17 problems in Table~\ref{table1}.
The remaining 12 problems are all
solvable\footnote{It is important here that jumps are restricted to lie within $T_5$.
If we use an infinite board, starting from a
triangular configuration of pegs (with one missing),
the SAX count is no longer a valid pagoda function and
all 17 feasible problems are solvable.}.
Schwartz and Ahlburg \cite{Schwartz} give an alternate
proof of some of the unsolvable cases.

The SAX count is also useful when solving the puzzle by hand.
It is useful to understand which jumps may result
in loss in the SAX count, so that these jumps can be avoided.
Jumps which can decrease the SAX count are as follows:
\begin{packed_enumerate}
\item Jumps ending in a corner (a1, a5 or e5) always lose 1.
\item Jumps beginning from the interior three holes (b3, b4, or c4) are particularly problematic.
Each such jump loses either 1 or 2.
It is a challenge to find \textbf{any} solution which includes
such a jump (it is possible for the c5-complement).
\item Edge-to-edge jumps, such as a3-c3, or a3-c5.
Edge-to-edge jumps do not always reduce the SAX count,
it depends on the number of pegs along the edges they connect.
These jumps are particularly important in the game since they are usually
the only way to reduce the number of interior pegs.
\item Jumps ending at one of the interior holes (b3, b4, or c4).
These jumps may or may not reduce the SAX count,
depending on the state of the rest of the board.
\end{packed_enumerate}

A player can avoid jumps of the first two types
by remembering the following rule of thumb:
\textit{avoid jumping into a corner or out of the interior}.
Many people stymied by this puzzle can find a solution
if they follow this rule of thumb.

\vskip 12pt

We can quantify how much this rule of thumb helps
by simulating players that select jumps at random.
At any board position, \textbf{Player A} counts the total number
of jumps available and selects one at random\footnote{There are other ways to
specify a ``random jump", and many give different behaviors.
For example, we could select a peg at random,
and then select a random jump using this peg
(repeating, if the selected peg has no jump).
This scheme gives a preference to jumps made by pegs
that have fewer jumps available.}.
\textbf{Player B} also selects a random jump,
but she follows the rule of thumb and does not consider
jumps into a corner or out of the interior.
Exceptions to the rule of thumb are made in the case of the first jump,
if no other option is possible (starting with a1 vacant, for example),
or on the last jump, if it ends with one peg.
Note that Player B may still make a jump which reduces the SAX count
to a level where a one peg solution can no longer be reached.

\textbf{Player C} calculates the SAX count after every potential jump,
and will not choose a jump that lowers the SAX count to less than
that of the finishing hole.
Since the player is aiming to finish with one peg anywhere,
Player C's best strategy is to keep the SAX count 
greater than or equal to zero, although the last jump
is allowed to violate this provided it ends with one peg.
In addition, for the b3 vacancy,
the SAX count begins at $-1$, so in this case
Player C considers only jumps that leave the SAX count unchanged.
In no case is Player B or C allowed to violate their rules,
unless there are 14 pegs on the board (first jump) or
2 pegs on the board (last jump).
If the number of pegs on the board is between 14 and 2,
and there is no jump that satisfies their rules, these
players have lost and must start over.

\begin{table}[htb]
\begin{center} 
\begin{tabular}{| c | c | c | c |}
\hline
$T_5$ & \multicolumn{3}{|c|}{Odds of finishing at a one peg position} \\
Starting & Player A & Player B & Player C \\
Vacancy & (any jump) & (rule of thumb) & (follows SAX count) \\
\hline
\hline
a1 & 1 in 146 & 1 in 47 & 1 in 7 \\ \hline
b3 & 1 in 579 & 1 in 142 & 1 in 19 \\ \hline
a4 & 1 in 291 & 1 in 47 & 1 in 7 \\ \hline
c5 & 1 in 141 & 1 in 13 & 1 in 7 \\ \hline
\end{tabular}
\caption{Odds of finishing with one peg for three players following different strategies
(the odds are not exact integers, they have been rounded to the nearest integer).} 
\label{table2}
\end{center} 
\end{table}

A computer can calculate the probability that these three players
finish at a one peg position.
This can be done by simulating a lot of games and accumulating
statistics, or we can calculate the exact probability that
each board position occurs, based on the game tree.
Table~\ref{table2} shows the results from an exact calculation.
The results show that a player can improve the odds of
finishing with one peg by at least a factor of 3 by using
the rule of thumb.
Much better is to calculate the SAX count exactly,
where an improvement by a factor of at least 20 is guaranteed.
Of course, human players do not select jumps at random,
and the result of human games would give even better odds than
those in Table~\ref{table2}.
Knowledge of the rule of thumb or (better) calculating the
full SAX count will clearly benefit human players as well.

It is interesting to note that players selecting jumps at random
have a much more difficult time with the standard
33-hole cross-shaped board.
Beginning from the usual position with one peg missing in the center,
Player A has only a 1 in 37 million chance of finishing with one peg
\cite{Bellweb}!
This is a consequence of the larger board size.

\section{Solving the Triangular Board \texorpdfstring{$T_n$}{Tn}}

 From any board position in the position class $\mbox{EMPTY}$,
it is impossible to play and finish with one peg.
In general, a board position in $\mbox{PEG}_i$
\textit{may or may not} be solvable to one peg.
For example the position with only holes a1 and a3 occupied
by pegs is in position class $\mbox{PEG}_1$, 
yet it cannot be reduced to a single peg.
However, in the case where $T_n$ is 
completely filled with pegs aside from a single missing peg,
we will now show that a board position in $\mbox{PEG}_i$
\textit{is always} solvable to one peg (provided $n\ge 4$).
This will supply the missing part of the proof of Theorem~\ref{th1}.

\subsection{Purges or Block Removals}

\begin{figure}[htb]
\centering
\epsfig{file=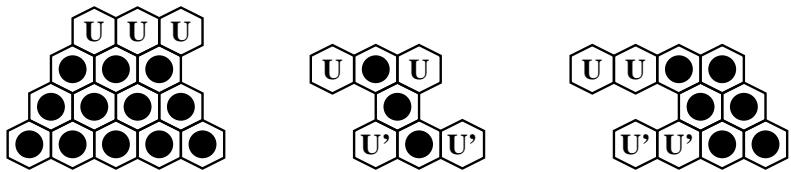}
\caption{Purges that work well together in triangular solitaire: trapezoid purge, 3-purge, 6-purge.}
\label{fig5}
\end{figure}

In square lattice solitaire it is useful to know block removals or
purges \cite{WinningWays, Beasley}---these are sequences
of jumps that remove a whole block of pegs,
leaving the rest of the board unchanged.
Figure~\ref{fig5} shows some useful purges on a triangular grid.
In each case, the pegs shown are removed by the purge,
while the holes labeled \textbf{U} must be \textbf{unlike},
in other words they cannot all be empty, or all filled by pegs.
The effect of the purge is to remove all the pegs shown
and restore the unlike holes, called the catalyst.
Holes labeled \textbf{U'} represent an alternative catalyst.
No jumps are allowed involving holes not shown in the figures
(and either \textbf{U} or \textbf{U'} must be chosen for the catalyst).

The 3- and 6-purges are well known from the square lattice case;
the reader can easily reconstruct the sequence of jumps
that performs the desired function.
The trapezoid purge is a new purge that works only on
a triangular grid, and the jumps to solve it are
more complicated and difficult to remember
(especially for all six combinations of unlike holes).
This purge is key to extending solutions on triangular
boards---we now go into some detail on it.

\begin{figure}[htb]
\centering
\epsfig{file=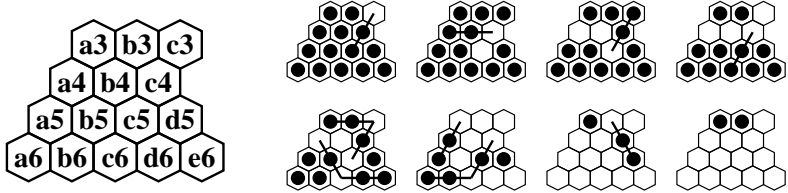}
\caption{A trapezoid purge starting from a3 and b3 filled, c3 empty.
Note that the jumps executed in reverse order solve the trapezoid purge with
a3 and b3 empty, c3 filled.}
\label{fig6}
\end{figure}

We can think of the trapezoid purge as
operating on the board $T_6$ with the six holes
a1, a2, b2, d4, e5, and f6 removed.
Note that the effect of this purge is to clear the bottom
three rows of this board, leaving the top row unchanged.
Figure~\ref{fig6} shows one way to execute the
purge which begins with a3 and b3 filled by pegs,
and c3 empty.
For solutions for other configurations of the unlike holes,
see Appendix~\ref{trap_sols}

The trapezoid purge can be extended to clear three consecutive
rows of any triangular board, provided the width of the top
row of the three is greater than 4.
We do this by stacking to the right of the trapezoid purge
at most one 3-purge and then as many 6-purges as necessary
(note that the 3-purge cannot be used at the edge of the board).
In order to ensure that the catalyst for these 3- and 6-purges is
available, it is useful to select trapezoid purge solutions
having the following property: at some time during their execution,
we have c4 empty, and b4 and e6 filled.
When we need the 3-purge, the fact that c4 is empty provides the
catalyst for this purge, and then the unlike pair for the 6-purges
appears on the bottom row.
If we only need 6-purges, the catalyst is provided by having c4 and
b4 unlike.
Note in Figure~\ref{fig6} that after the first move
the board has c4 empty and b4 and e6 filled.

\begin{figure}[htb]
\centering
\epsfig{file=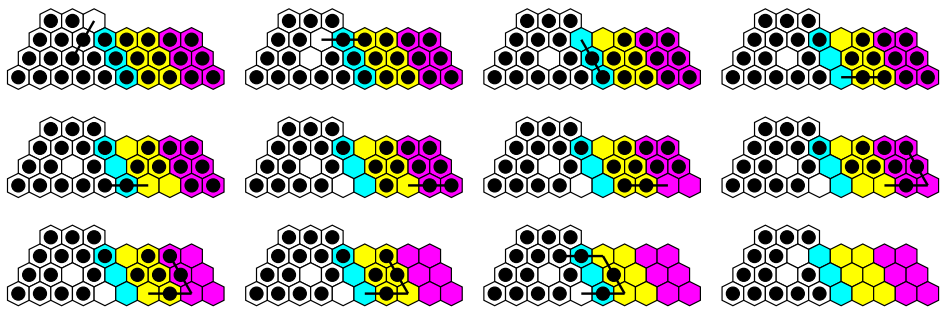}
\caption{Clearing the bottom three rows of $T_{10}$ using
a trapezoid, 3-purge and two 6-purges (11 moves, 16 jumps,
the finish of the trapezoid purge is not shown).}
\label{fig7}
\end{figure}

An example using the bottom three rows of $T_{10}$
is shown in Figure~\ref{fig7},
starting from the same trapezoid purge catalyst.
Note that the final board position in Figure~\ref{fig7}
is the same as that in Figure~\ref{fig6} after the first jump,
so the final jumps of the trapezoid purge are identical and
therefore not shown at the end of Figure~\ref{fig7}.
In Figure~\ref{fig7}, the first jump comes from the trapezoid purge,
while jumps 2,3 and 14 come from the 3-purge\footnote{Here we count individual jumps,
rather than moves, because a move may contain jumps from different purges.}.
The left 6-purge is done on jumps 4, 5, 12, 13, 15, and 16, and the
right 6-purge is jumps 6--11.
The exact sequence of jumps may seem hard to figure out,
but in the next section we will give a simple
algorithm to determine the sequencing.

This technique of clearing three consecutive rows
gives us an inductive technique for extending
solutions on triangular boards.
For example, suppose we have a solution on $T_7$,
beginning with one peg missing and ending with one peg.
This can be extended to a solution on $T_{10}$ by clearing
the bottom three rows as in Figure~\ref{fig7}
(although not necessarily using the
Figure~\ref{fig7} trapezoid purge catalyst).
Note that some catalyst for the trapezoid purge
\textit{must always be present} at some time during
the $T_7$ solution.
If the starting vacancy is one of the three catalyst holes,
then the catalyst is present at the start.
Otherwise, it is impossible for a single jump to remove all three
pegs in the catalyst area.
The fact that the catalyst must always be present is a critical
feature of this trapezoid purge,
and is not the case for many other purges involving
only two unlike holes.

\subsection{A Simple Algorithm for Solving \texorpdfstring{$T_n$}{Tn}}
\label{finish_th1}

Let us consider the general question posed in Theorem~\ref{th1}:
when can a feasible vacancy on $T_n$ be solved to one peg?
Note here that we are not free to select the location of the
final peg---this more specific case will be handled in the
next section.
The strategy is to use solutions on $T_4$, $T_5$ and $T_6$ to
inductively define solutions for all larger boards.

First, we note that many starting vacancies are equivalent
by rotation and/or reflection of the board (see Appendix~\ref{rotref}).
Therefore we need only a few solutions on $T_4$, $T_5$ and $T_6$
to be able to solve all vacancies on them.
In Appendix~\ref{Tn_sols}, we give solutions that cover all
problems on $T_4$, $T_5$ and $T_6$.

We now consider any feasible starting vacancy
on $T_n$ with $n>6$.
First, we choose a sub-board
among $T_4$, $T_5$ and $T_6$ that has the same remainder when divided by 3.
We now place this sub-board inside the larger board so that both:
\begin{packed_enumerate}
\item The sub-board encloses the starting vacancy.
\item The number of holes between the edges of the sub-board and
larger board is a multiple of 3 (in all three directions).
\end{packed_enumerate}

\begin{figure}[htb]
\centering
\epsfig{file=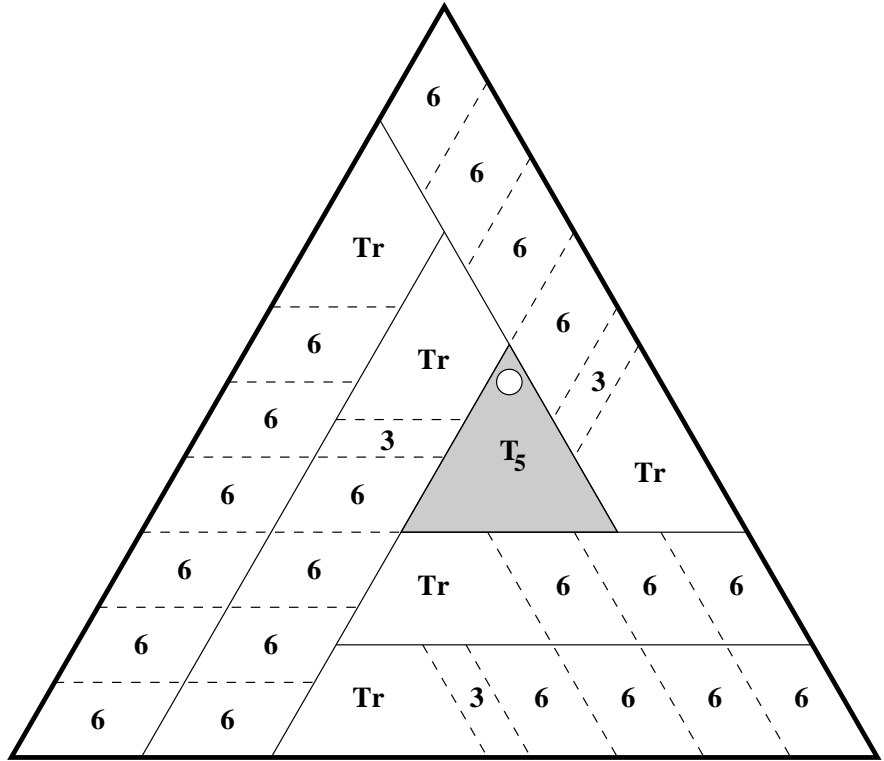}
\caption{Solving the g10 vacancy on $T_{20}$ using a solution on the $T_5$ sub-board and purges.
Regions marked by ``Tr", ``3" and ``6" are cleared by trapezoid, 3, or 6-purges, respectively.}
\label{fig8}
\end{figure}

An example of the decomposition of the $T_{20}$ board for a
$\mbox{g10}=(6,9)$ starting vacancy is shown in Figure~\ref{fig8}.
Since $20\equiv 2 \pmod 3$,
we select $T_5$ as the sub-board, and we will use the
solution to the a1 vacancy on this board
(note that alternatively we could have placed the top
corner of the sub-board at d7 or g7).
The remaining portions of the board are cleared by appropriate purges,
as mapped out in Figure~\ref{fig8}.
We have been careful in designing our trapezoid purges to ensure that
the required catalyst for each purge
is available at some time during the solution.

Working by hand, it is not trivial to find a sequence of jumps
for the solution diagrammed in Figure~\ref{fig8}.
As each set of three parallel rows is cleared,
the purges become interleaved in the sense that the next one
begins before the last one finishes (as in Figure~\ref{fig7}),
and this can be difficult to keep track of.
The algorithm is well suited for a computer, however,
and runs extremely quickly.
We have programmed up an online triangular game on the web \cite{Bellgame}
that can solve any feasible vacancy on $T_n$ for $n\ge 4$
(although due to display limitations it will only go up to $T_{24}$).

Here is the algorithm used in my program \cite{Bellgame}.
First, find the size of the sub-board and its location.
Then by rotation and/or reflection of a known solution,
we obtain a solution to the problem on the sub-board.
We then determine a list of the purges to clear the rest
of the board (for the problem of Figure~\ref{fig8},
this list would have 29 purges).
Associated with each purge is a set of unlike catalyst holes,
and a counter that indicates how many jumps in that purge have been executed
(and if it is finished).
We now initialize the board at the starting position and execute
the following algorithm to determine the sequence of jumps:

\begin{packed_enumerate}
\item Go through the purges with 0 jumps executed.
If the catalyst for this purge is present,
execute the first jump in this purge and return to Step 1.

\item Go through the purges that have been started but
are not finished in the reverse order that they were started in.
If the next jump in a purge is possible
(pegs in the correct configuration),
execute that jump and return to Step 1.
Otherwise check the next purge.

\item If the board now contains only one peg, stop,
this is the final board position.
Otherwise, execute the next jump in the sub-board solution.
Return to Step 1.
\end{packed_enumerate}

The reader can check that this algorithm gives
the same sequence of jumps shown in Figure~\ref{fig7}.
The reader is also urged to watch this solution technique on the web
version of the puzzle \cite{Bellgame}.
Note that the algorithm \cite{Bellgame} uses a wider trapezoid purge
(one hole wider) to extend solutions on $T_5$, $T_6$ and $T_7$.

\subsection{A More General Algorithm for Solving \texorpdfstring{$T_n$}{Tn}}
\label{finish_th2}

In this section we will prove that as long as the board size $n\ge 6$,
it is possible to play between any feasible pair of starting and ending holes.
This will supply the missing half of the proof of Theorem~\ref{th2}.
Table~\ref{table3} lists the number of distinct feasible problems on each board,
the formula for general $n$ is most easily derived using
Burnside's Lemma.

The proof is inductive on $n$, and the first step is to verify that
all problems can be solved for $n=6$, $7$, or $8$.
This is non-trivial due to the large number of problems,
particularly for $n=8$ where there are $80$ cases.
One useful trick is that if we have a solution from
$(x_s,y_s)$ to $(x_f,y_f)$, then by playing the jumps in
the reverse order we obtain a solution from $(x_f,y_f)$
to $(x_s,y_s)$.
Many cases can also be covered by extending solutions on $T_5$.

\begin{table}[htb]
\begin{center} 
\begin{tabular}{| c | c | c | c |}
\hline
 & & \multicolumn{2}{|c|}{Number of Feasible Pairs}\\
Board Side & Board Size & Distinct & Solvable\\
($n$) & $T(n)$ & \seqnum{A130515} \cite{OEIS} & \seqnum{A130516} \cite{OEIS}\\
\hline
\hline
2 & 3 & 1 & 0\\ \hline
3 & 6 & 4 & 0\\ \hline
4 & 10 & 3 & 1\\ \hline
5 & 15 & 17 & 12\\ \hline
6 & 21 & 29 & 29\\ \hline
7 & 28 & 27 & 27\\ \hline
8 & 36 & 80 & 80\\ \hline
9 & 45 & 125 & 125\\ \hline
10 & 55 & 108 & 108\\ \hline
11 & 66 & 260 & 260\\ \hline
12 & 78 & 356 & 356\\ \hline
\multicolumn{2}{|c|}{$n\equiv 1 \pmod 3$} & \multicolumn{2}{|c|}{$(T(n)-1)^2/27$}\\ \hline
\multicolumn{2}{|c|}{$n\not\equiv 1 \pmod 3$ and $n$ even} & \multicolumn{2}{|c|}{$(4T(n)^2+9n^2)/72$}\\ \hline
\multicolumn{2}{|c|}{$n\not\equiv 1 \pmod 3$ and $n$ odd} & \multicolumn{2}{|c|}{$(4T(n)^2+9(n+1)^2)/72$}\\ \hline
\end{tabular}
\caption{The number of distinct feasible pairs starting with one peg missing and finishing with one peg.
Theorem~\ref{th2} states that for $n\ge 6$ the rightmost two columns are equal.} 
\label{table3}
\end{center} 
\end{table}

We now consider the inductive step.
We want to show that any feasible problem on $T_n$ is solvable if we know
that every feasible problem on $T_{n-3}$ is solvable.
If $(x_s,y_s)$ and $(x_f,y_f)$ are close enough together,
we can enclose them inside $T_{n-3}$ where one of the three corners of
$T_{n-3}$ coincides with a corner of $T_n$.
We can then play our solution on $T_{n-3}$, and clear
the remainder of the board using trapezoid, 3- and 6-purges,
exactly as in the previous section.

The remaining case, therefore, occurs when $(x_s,y_s)$ and $(x_f,y_f)$
are far enough apart that they \textit{cannot} be contained inside $T_{n-3}$.
Let us rotate and/or reflect the board so that the starting vacancy
$(x_s,y_s)$ is in the last three rows of $T_n$, and the finishing
hole $(x_f,y_f)$ is not in the last three rows.

\begin{figure}[htb]
\centering
\epsfig{file=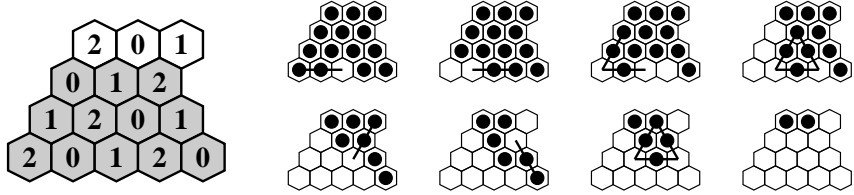}
\caption{A specific purge to clear the last three rows.  The numbers on the left show the position class of each hole.}
\label{fig9}
\end{figure}

What we need now is a specific type of purge that can empty the last three
rows of $T_n$, and leave all rows above filled except for a
single vacancy
(which must be in the same position class as the starting vacancy).
This is easy to accomplish using a type of trapezoid purge as
diagrammed in Figure~\ref{fig9}a.
The starting board position has every hole filled by a peg except for one
vacancy in the shaded region.
The target board position has the bottom three rows empty, and
a vacancy in the top row only at the hole of the same parity type as
the starting vacancy.
Figure~\ref{fig9}b shows an example of such a solution.
The solution to this purge for other starting vacancies is left as an
exercise for the reader.

To clear the remainder of the bottom three rows,
we insert 3- and 6-purges to the right at appropriate times.
If the starting vacancy is farther to the right, we translate some purge of
Figure~\ref{fig9} to the right, or reflect it, with the left and right portions
cleared again by combinations of 3- and 6-purges.
As before we select our trapezoid purge solution with
certain properties\footnote{To insert purges on the right,
it suffices to have at some time c4 empty and b4 and e6 filled,
\textit{or} e6 empty and c4 and d6 filled.
For the purges on the left, we need a4 empty and b4 and a6 filled,
\textit{or} a6 empty and a4 and b6 filled.}
to ensure the catalyst
for the 3- and 6-purges is present at some time.
The board above the bottom three rows is now solved using a solution on $T_{n-3}$.

\section{Short Solutions}

The previous section showed how to solve any feasible
combination of starting vacancy and finishing hole
for board side $n\ge 6$.
A much more difficult task is to find the \textit{shortest}
solution to such a problem.
Here by shortest we mean a solution with the minimum number of moves,
where a move is one or more consecutive jumps by the same peg.

\subsection{Bounds}

Let $S(n)$ denote the length of the shortest solution
(in moves) to any problem on $T_n$ starting with
one peg missing and ending with one peg
($S(n)$ is \seqnum{A127500} in OEIS \cite{OEIS}).
A trivial upper bound on $S(n)$ is the number of jumps
in any solution, which is the size of the board minus two,
$T(n)-2=n(n+1)/2-2$.

\begin{figure}[htb]
\centering
\epsfig{file=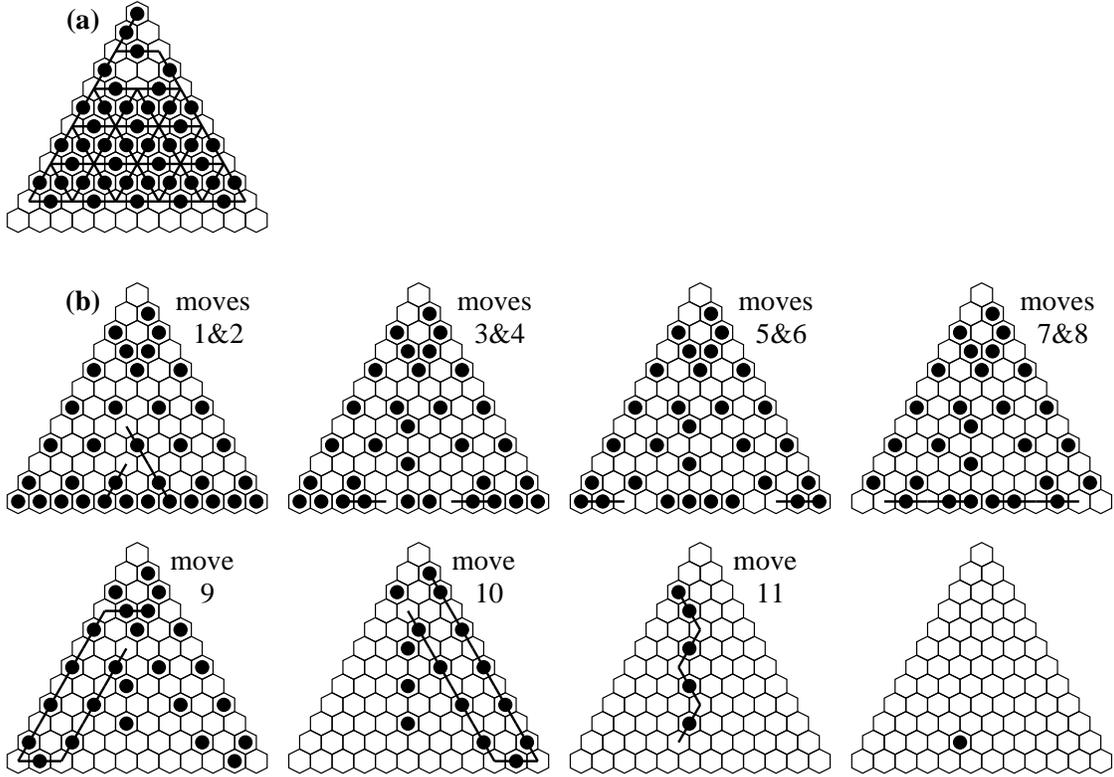}
\caption{Building a short solution on $T_{12}$:
(a) a final move which removes 42 pegs
(b) reducing the complement of (a) to one peg.}
\label{fig10}
\end{figure}

We can obtain a better upper bound on $S(n)$ by finding
a solution with a final move with as many jumps as possible.
Figure~\ref{fig10}a shows such a move on $T_{12}$,
this move begins from a1, removes 42 pegs, and finishes
with one peg at a3.
But can this be the final move to a peg solitaire problem?
In other words, is it possible to reach the board position in
Figure~\ref{fig10}a starting from a full board with
one peg missing?

To answer this question, we can use an elegant technique called
the ``time reversal trick" in Winning Ways \cite[p. 817--8]{WinningWays}.
Figure~\ref{fig10}b shows how to start from the 
complement of the board position in Figure~\ref{fig10}a,
and reduce the board to one peg.
If we take the complement of the final board position in
Figure~\ref{fig10}b, and execute the jumps in reverse order,
we will reproduce the board position in Figure~\ref{fig10}a.
When the jumps are played in reverse order, they are generally
broken into separate moves, so this solution contains
34 single jump moves,
plus the long final move, for a total of 35~moves.
However, by cleverly reordering the first 34~jumps
we can create a solution with only 29~moves \cite{GPJ04}.

The remarkable thing about the solution in Figure~\ref{fig10}
is that it can be extended to
any \textit{even} board size $n\ge 12$.
For $12<n<24$,
we simply extend the final sweep pattern to fill the bottom
of the board.
To solve the problem analogous to Figure~\ref{fig10}b,
we extend moves 1, 2, 9, and 10, and add extra moves
after move 11 to finish with one peg.
For $n>22$ we need to show how to reduce
the lower, central portion of the board to one peg.
Bell \cite{GPJ04} gives the details of this inductive step.
The net result is an upper bound on $S(n)$,
\begin{equation}
S(n) \le \frac{1}{8}n^2+\frac{7}{6}n-3.
\label{UBound}
\end{equation}
The upper bound (\ref{UBound}) is valid when $n$ is a multiple of 12.
For any even $n\ge 12$ the construction
analogous to that in Figure~\ref{fig10}
gives a bound with the same leading order term,
but different linear and constant terms.

\begin{figure}[htb]
\centering
\epsfig{file=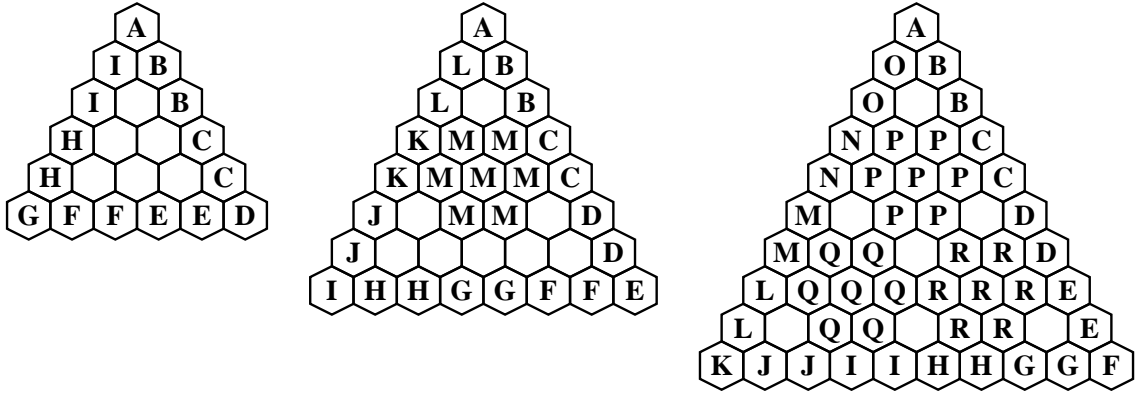}
\caption{Merson regions on $T_6$, $T_8$ and $T_{10}$, $r=9$, $13$, and $18$, respectively.}
\label{fig11}
\end{figure}

A lower bound on $S(n)$ can be obtained from the following argument:
consider the board $T_n$ divided into $r$
``Merson regions\footnote{Named after Robin Merson who
first used this concept in 1962 on the $6\times 6$
square board \cite[p. 203]{Beasley}.}"
as in Figure~\ref{fig11}.
The shape of a region is chosen such that when it is entirely filled with pegs,
there is no way to remove a peg in the region without a move that originates
inside the region.
Each of the three corners is a region,
as well as any pair of consecutive holes along the edge.
In the interior of the board the regions must be
large 7-hole hexagons.

Any region that starts out filled must have at least one
move starting from inside it.
Since the starting position has every hole filled by a peg except one,
all regions start filled except possibly the region that
contains the starting hole.
If the board can be divided into $r$ regions,
then no solution beginning with one peg missing and ending at one peg
can have fewer than $r-1$ moves.
\begin{equation}
T \left( \left \lfloor \frac{n-4}{3} \right \rfloor \right)
+ \left \lfloor \frac{3n-2}{2} \right \rfloor
\le S(n),
\label{LBound}
\end{equation}
where $T(n)=n(n+1)/2$ is a triangular number.
For $n$ even,
the lower bound in (\ref{LBound}) is $r-1$.
For $n$ odd, there will always be a gap between the edge regions,
and often we can choose the regions so the starting vacancy
is not in any region.
Even when it is not,
by considering the first few moves, the lower
bound can be taken as $r$ for $n$ odd,
as given in (\ref{LBound}).

When $n$ is a multiple of 12,
the upper bound (\ref{UBound}) and lower bound (\ref{LBound}) reduce to
\begin{equation}
\frac{1}{18}n^2+n \le S(n) \le  \frac{1}{8}n^2+\frac{7}{6}n-3.
\label{ULBounds}
\end{equation}
For all even $n\ge 12$, the upper and lower bounds
have the same leading order behavior as given in (\ref{ULBounds}).

\subsection{Computational Search Techniques}

We now turn to determining the value of $S(n)$ by computer search.
An exhaustive search for short solutions is difficult
beyond $T_7$.
A more efficient search technique is
``breadth-first iterative deepening A*"
as defined by Zhou and Hansen \cite{ZhouAI}
and applied to peg solitaire by Bell \cite{BellDiag}.
This algorithm can find $S(n)$ up to $n=10$,
with results shown in Table~\ref{table4}.
Solutions of length $S(n)$ can be found in Appendix~\ref{Tn_sols}.

\begin{table}[htb]
\begin{center} 
\begin{tabular}{| c | c | c | c |}
\hline
Board Side & Board Size & Lower & Shortest\\
($n$) & $T(n)$ & Bound (\ref{LBound}) & Solution $S(n)$\\
 & & & \seqnum{A127500} \cite{OEIS}\\
\hline
\hline
4 & 10 & 5 & 5\\ \hline
5 & 15 & 6 & 9\\ \hline
6 & 21 & 8 & 9\\ \hline
7 & 28 & 10 & 12\\ \hline
8 & 36 & 12 & 13\\ \hline
9 & 45 & 13 & 16\\ \hline
10 & 55 & 17 & 18\\ \hline
11 & 66 & 18 & $19 \le S(11) \le 28$\\ \hline
12 & 78 & 20 & $21 \le S(12) \le 29$\\ \hline
\end{tabular}
\caption{The minimum number of moves needed to solve problems on $T_n$.} 
\label{table4}
\end{center} 
\end{table}

Determining the value of $S(n)$ computationally
involves several steps.
For example, for $n=10$, we 
first apply the search algorithm to look for solutions
of length 17 (or less) at all geometrically
distinct starting vacancies from among
$108$ feasible pairs\footnote{For $T_{10}$, there are twelve geometrically
distinct starting vacancies.
However, $10\equiv 1 \pmod 3$,
so an application of Theorem~\ref{th2} finds that six of these cannot be reduced to one peg.
So for $n=10$ we need only check
the starting vacancies: a2, a3, a5, b4, b5 and c6.}.
If all searches finish with no solution found,
we know $S(10)>17$.
The next step is to run the search algorithm to find
a solution of length 18.
This is the most time consuming step.

For $n=11$, our search algorithm can determine
that no solution of length 18 exists.
Searching for a solution of length 19 is, however,
too difficult, and the algorithm runs for several
days before running out of disk space.
The shortest known solution for $n=11$ has 28 moves \cite{Bellweb}.

\section{Conclusions}

In this paper, we have considered the special problem where
the board begins from a position with one peg missing,
and finishes with one peg.
We have found necessary and sufficient conditions for such problems
to be solvable on any size triangular board.
Moreover, for all solvable problems we have given a fast
solution algorithm which does not rely on exhaustive search,
but builds a solution in an inductive manner using
solutions on smaller boards (which are pre-computed
for the smallest size boards).
These solution extension techniques can also be applied to
other board shapes besides triangular, such as rhombus \cite{GPJ05}
and hexagonal.

Similar ideas for extending solutions can also
be applied in square lattice solitaire.
A modified trapezoid purge is needed for this case,
and one possibility is the 15-hole board formed by
a $4\times 4$ square board with the upper right
corner removed, with the three
holes in the top row forming the unlike catalyst.
As in Section~\ref{finish_th1}, solutions on rectangular boards
can be extended using this new purge together with
the usual 3- and 6-purges.

A more general problem is to determine if an \textit{arbitrary}
configuration of pegs can be reduced to a single peg.
It has been proved (for the case of an $n\times n$ square board)
that this problem is NP-complete \cite{Uehara}.
 From a practical standpoint this means we cannot hope to find
a fast (polynomial speed) algorithm for solving this general problem.
What we have presented here is a small subset of problems that can be
solved much more quickly.

Finally, we have considered the problem of finding the shortest solution
(in moves) to any problem beginning with one peg missing and
finishing with one peg.
This is a much more difficult task than finding any solution,
and we have shown how it can be solved (up to $T_{10}$)
using computational search techniques.

\section*{Acknowledgments}
I thank Pablo Guerrero-Garc{\'i}a
of the University of M{\'a}laga (Spain)
and the anonymous referee,
for many useful comments which improved this paper.

\appendix
\section{Skew Coordinate Transformations}
\label{rotref}

Many starting board positions are equivalent by reflection and/or rotation of the board.
For the triangular board $T_n$, we consider a reflection $f$ about the vertical axis,
and a rotation $r$ counter-clockwise by $120^\circ$.
These transformations are useful inside programs to convert solutions.
In skew coordinates these two coordinate transformations have a simple form,
given by
\begin{eqnarray*}
f(x,y) & = & (y-x, y), \\
r(x,y) & = & (y-x, n-1-x).
\end{eqnarray*}

Note that $f^2=i$ and $r^3=i$, where $i$ is the identity transformation.
There are six different transformations of the board,
given by the generators $\{i, r, r^2, f, rf, r^2f \}$.
For example, on $T_5$ one feasible pair is $(a2, b4)$, in
skew coordinates this pair is $((0,1), (1,3))$.
By applying a rotation, we obtain an equivalent feasible pair,
$(r(0,1),r(1,3))=((1,4),(2,3))=(b5,c4)$.
The six transformations above generate the equivalent set of
six feasible pairs (converted to alphanumeric notation):
$\{(a2,b4), (b5,c4), (d4,b3), (b2,c4), (a4,b3), (d5,b4)\}$.

\section{Solutions}

\subsection{Trapezoid Purges}
\label{trap_sols}

Trapezoid purge solutions for Section~~\ref{finish_th1},
using the coordinate system in Figure~\ref{fig6}:\newline
a3 empty, b3 and c3 filled:
a5-a3, c6-a4, a3-a5, a6-a4, c3-a3-a5, c5-c3, e6-c6-a6-a4-c4, d5-b3;
a3 and c3 filled, b3 empty:
d5-b3, c6-c4, c3-c5, e6-c6-c4, a3-c3-c5, a5-a3, a6-c6-a4-c4, c5-c3;
a3 and b3 filled, c3 empty (Figure~\ref{fig6}):
c5-c3, a4-c4, c3-c5, c6-c4, a3-c3-c5, e6-c6-a4, a5-a3, a6-c6-c4, d5-b3.

For the other three cases when the catalyst is
in the complementary configuration,
play the jumps in the appropriate solution above in reverse order.
For example, if a3 is filled, and b3 and c3 are empty, play:
d5-b3, a4-c4, a6-a4, c6-a6, e6-c6, c5-c3, $\ldots$.

\subsection{\texorpdfstring{$T_n$}{Tn} Solutions}
\label{Tn_sols}

For $T_4$, $T_5$ and $T_6$, we present a set of solutions,
which (appropriately rotated and/or reflected)
can reduce any feasible starting vacancy to one peg.
For $7\le n \le 10$ we give one solution on $T_n$
with the minimum number of moves $S(n)$.
Most of these solutions were found by computer,
for more information see my web site \cite{Bellweb}.
For diagrams of the larger solutions see
\href{http://www.geocities.com/gibell.geo/pegsolitaire/LargeTriangular/}
{\tt http://www.geocities.com/gibell.geo/pegsolitaire/LargeTriangular/}

\vskip 12pt 
\noindent\textbf{$T_4$ Solution}:\newline
Vacate a2: a4-a2, a1-a3, c4-a4-a2, c3-a3-a1-c3, d4-b2 (5 moves).

\begin{figure}[htb]
\centering
\epsfig{file=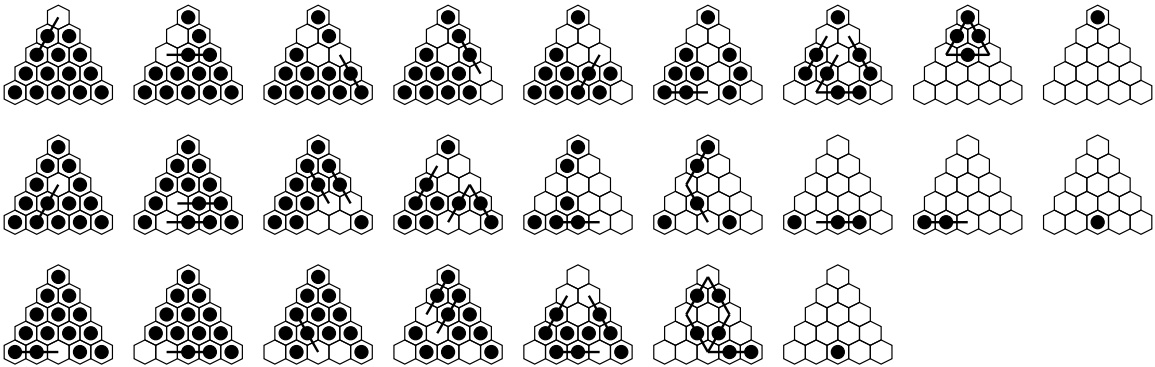}
\caption{Solutions on $T_5$.  The solution in the top row defines \seqnum{A120422} (after conversion to an integer hole notation).}
\label{figT5}
\end{figure}

\vskip 12pt\noindent 
\textbf{$T_5$ Solutions}:\newline
Vacate a1: a3-a1, c3-a3, e5-c3, b2-d4, c5-c3, a5-c5, d5-b5-b3, d4-b2, a4-a2, a1-a3-c3-a1 (10~moves, Figure~\ref{figT5}a, \seqnum{A120422});
vacate a4: a2-a4, then continue as previous solution;
vacate b3: b5-b3, d4-b4, d5-b5, b2-d4, a2-c4, a4-a2, e5-c3-c5, b5-d5, a1-a3-c5, d5-b5, a5-c5 (11~moves, Figure~\ref{figT5}b);
vacate c5: a5-c5, d5-b5, a3-c5, a1-a3, b2-b4, d4-b2, a4-a2, b5-d5, e5-c5-c3-a1-a3-c5 (9~moves, Figure~\ref{figT5}c).

\begin{figure}[htb]
\centering
\epsfig{file=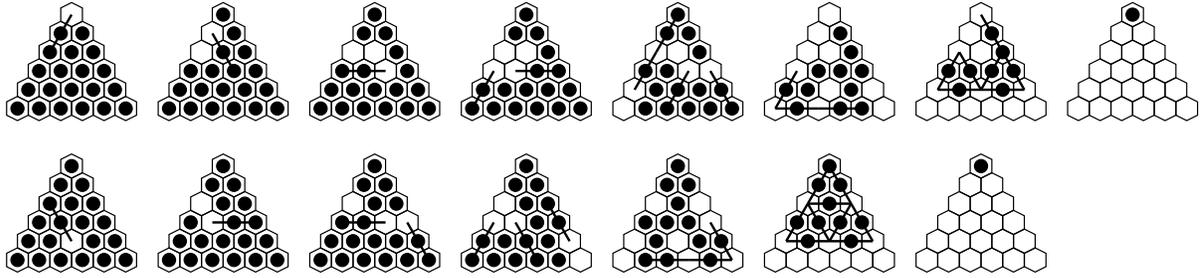}
\caption{Solutions on $T_6$.}
\label{figT6}
\end{figure}

\vskip 12pt\noindent 
\textbf{$T_6$ Solutions}:\newline
Vacate a1: a3-a1, c4-a2, a4-c4, d4-b4, a6-a4, a1-a3-a5, c6-c4, f6-d4, e6-c6-a6-a4, c3-e5-c5-a5-a3-c5-c3-a1 (10~moves, Figure~\ref{figT6}a);
vacate a4: a6-a4, a3-a5, a1-a3, c4-a2-a4-c4, d4-b4, c6-c4, e6-c6-a6-a4, f6-d4, c3-e5-c5-a5-a3-c5-c3-a1 (9~moves);
vacate b3: d5-b3, c6-c4, c3-c5, a6-c6, d6-b6, f6-d6, a4-c4, a2-a4-a6-c6-e6, a1-c3, d4-f6-d6-b4-b2-d4-b4-b6 (10~moves);
vacate c5: a3-c5, d4-b4, a4-c4, f6-d4, a6-a4, c3-e5, d6-b4, b6-d6-f6-d4, a1-a3-a5-c5-e5-c3-c5-a3-c3-a1 (9~moves, Figure~\ref{figT6}b);
vacate b6: d6-b6, a6-c6, f6-d6-b6, c4-e6, a4-a6-c6-c4, c3-c5, a2-a4-c4, a1-c3, d4-b4-b2-d4-f6-d6-b4-b6 (9~moves).

\vskip 12pt\noindent
\textbf{$T_7$ Solution}:\newline
Vacate c3: a1-c3, d4-b2, f6-d4, a3-c3-e5, d6-d4-f6, b4-d6, a5-c5, f7-d5-b5, d7-f7, g7-e7, b7-d7-f7, a7-a5-c7-c5-a5-a3-a1-c3-c5-e7-g7-e5 (12~moves).

\begin{figure}[htb]
\centering
\epsfig{file=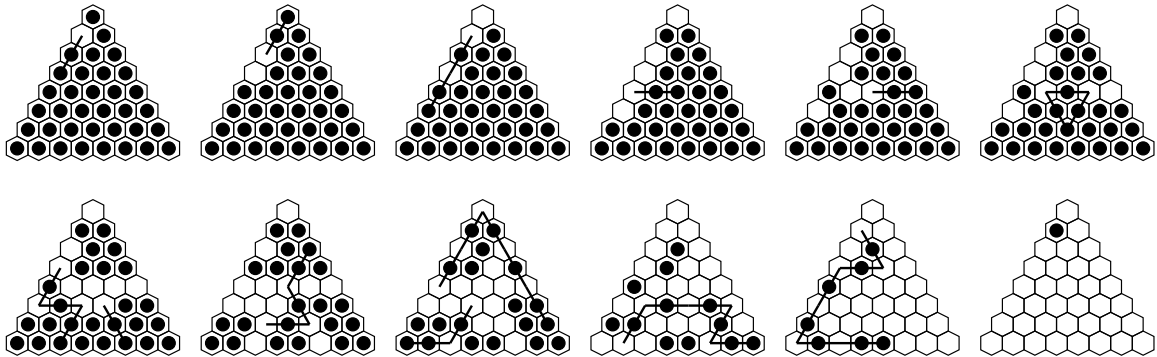}
\caption{A 13-move solution on $T_8$.}
\label{figT8}
\end{figure}

\vskip 12pt\noindent 
\textbf{$T_8$ Solution} (the only complement problem solvable in 13~moves):\newline
Vacate a2: a4-a2, a1-a3, a6-a4-a2, c5-a5, e5-c5, d7-d5-b5-d7, c8-c6-a6-a4, f8-d6,
c3-c5-e7-c7, a8-c8-c6, g7-e5-c3-a1-a3-a5, h8-f8-f6-d6-b6-b8, e8-c8-a8-a6-a4-c4-a2 (13~moves, Figure~\ref{figT8}).

\vskip 12pt\noindent 
\textbf{$T_9$ Solution}:\newline
Vacate a2: a4-a2, a6-a4, c5-a3-a5, e7-c5, g9-e7, d4-d6-f8, i9-g9-e7,
f6-d4-b4-d6-f8, c7-c5, a1-a3, h8-f6, e9-e7, c9-c7, a9-c9-e9-g9-g7-e5,
b2-d4-f6-d6-b4, a8-c8-e8-g8-e6-e8-c6-a4-c4-a2-a4-a6-c6-c8-a6-a8 (16~moves).

\vskip 12pt\noindent 
\textbf{$T_{10}$ Solution}:\newline
Vacate a3: a1-a3, a4-a2, a6-a4, a8-a6, c3-a1-a3-a5-a7, c5-a3, e5-c3-c5-a5, g7-e5-c5,
f8-f6, f10-f8, d7-d5-b5-d7-f9-f7-d5, c8-c6-a6-a8-c8-e8-e6, d10-b8-b6, b10-d10-f10-d8-d10,
i9-g7-e5-e7, h10-h8-f8, j10-h10-f10, a10-a8-c10-e10-g10-g8-e8-e6-c4-a2-a4-a6-c6 (18~moves).

\bigskip
\hrule
\bigskip

\noindent 2000 {\it Mathematics Subject Classifications}:
Primary 00A08; Secondary 97A20.

\noindent \emph{Keywords and phrases:} triangular peg solitaire,
marble solitaire puzzle, peg jumping game.

\bigskip
\hrule
\bigskip

\noindent (Concerned with sequences
\seqnum{A120422}, \seqnum{A127500}, \seqnum{A130515}, and \seqnum{A130516}.)

\bigskip
\hrule
\bigskip

\vspace*{+.1in}
\noindent
Received August 1 2008;
revised version received November 12 2008. 
Published in {\it Journal of Integer Sequences}, November 16 2008.

\bigskip
\hrule
\bigskip


\noindent\textbf{\Large Addendum}\newline
Shortly after this paper was published in
\htmladdnormallink{JIS}{http://www.math.uwaterloo.ca/JIS/},
I received an email from John Beasley with
many useful comments.
Further correspondence
simplified several of the arguments in this paper.
I have condensed our discussion into the comments below.

In order to interpret John's comments, one needs
to understand the ``exit theorems", a concept
introduced in his book \cite[Chapter 7]{Beasley}.
Consider a region (any fixed subset of holes).
We call a jump an \textbf{exit} from this region
if it removes a peg from the region and
ends outside the region.
If a region is not empty, then clearly to remove all pegs
from this region we must have at least one exit
(this will be the jump which leaves the region empty).
Similarly, if a region consists of three or more holes
and begins full and finishes empty,
we must have at least two exits
(one to remove the first peg and one to remove the last peg).
These seemingly obvious statements are called
the \textbf{exit theorems} \cite[p. 117]{Beasley}
and prove surprisingly powerful.

\bigskip
\hrule
\bigskip

\noindent{\large Comments from John Beasley, December 2008:}\newline

\noindent\textbf{Section 2.2, $T_4$ problem ``vacate a2, finish at c4"}\newline
Consider the set of holes $\{a2, a4, c4\}$.
This is a closed set, in the sense that no peg can jump into it from outside,
so the number of pegs in these holes can never increase.
We start with two and we must end with one,
so we can only lose one in the course of the solution.

Now consider the region $\{a1, b2, c3, d4\}$.
This starts full and ends empty, so we need two exits from it.
However, by your own parity argument,
the move across b3 must be a2-c4 or c4-a2,
so neither b2-b4 nor c3-a3 is available,
and the only other possibilities,
a1-a3 and d4-b4, each remove a peg from the set $\{a2, a4, c4\}$.
Hence the problem is unsolvable.

\newpage
\noindent\textbf{Section 2.3, the SAX count on $T_5$}\newline
The SAX count can be derived from an entirely different perspective
by counting exits to six carefully chosen regions.
Consider the six regions consisting of the three corners
$\{a1\}$, $\{a5\}$, $\{e5\}$ plus
the ``edge regions" $\{a2,a3,a4\}$, $\{b2, c3, d4\}$ and $\{b5, c5, d5\}$.
We call the pegs in
$\{a2$, $b2$, $b3$, $a4$, $b4$, $c4$, $d4$, $b5$, $d5\}$
the ``fodder pegs"; for any particular board position,
let $F$ be the number of fodder pegs.
The reason for this name is that
\textit{any exit from these six regions
necessarily consumes a fodder peg}.
Two important details are that
\textit{the number of fodder pegs can never increase} (they form a closed set),
and an exit from one region cannot also be an exit from
a different region\footnote{This might not be the case if, for example,
a region lies in the interior of the board.}.
This immediately provides a simple proof that
the b3-complement is unsolvable---for we need
9 exits from the six regions,
yet we have only 8 fodder pegs at the
start and must finish with 1 fodder peg.

In fact, from any board position we can do an accounting of
the remaining fodder pegs $F$ and subtract the number
of remaining exits,
I call this the ``$F-E$ count".
Specifically, to calculate $E$ we take the number of corners occupied,
and for each edge region, add +2 if this region is full,
+1 if the region contains 1 or 2 pegs, and 0 if it is empty.

Because each exit consumes a fodder peg,
and the number of fodder pegs can never increase,
this proves that the $F-E$ count can not increase.
In fact, it is not hard to show that the $F-E$ count is
\textit{identical} to the SAX count!

This embeds the mysterious SAX count in a more
general framework of regions and exits which is
applicable to any board.
The proof that it can never increase is
automatic and does not require a case-by-case
analysis as for the SAX count \cite{Hentzel2}.
In addition, remarks 3 and 4 lower down on page 7 (of this paper)
can be rephrased in terms of exits and the fullness or emptiness
of the three-hole edge region affected,
when they are perhaps more readily comprehended.
 
\bigskip\noindent\textbf{Section 3.1, the trapezoid purge}\newline
I think you are doing yourself a slight injustice in calling this
``more difficult to remember".  Consider the following:
Since UUU are unlike, at least one of them must be empty.
This gives us a catalyst for the 3-purge a4/b4/c4.
 
By the same argument, at least one of UUU must be full.
Suppose first that a3 is full;
then we play c6-c4, e6-c6, b6-d6, c4-e6-c6, reducing to a
\textbf{hollow V} a6/a5/b5/c6 hinged on a4, and the combination
``a3 full, a4 empty" gives us a catalyst to purge it
(a6-a4, a3-a5, c6-a4, a5-a3).
If instead, b3 or c3 is full,
we play to leave the hollow V on c6/c5/d5/e6.
 
I think this combination ``purge a4/b4/c4,
then play to leave the appropriate hollow V"
gives a simple and easily remembered solution to
what is really rather an ingenious removal.

\end{document}